\documentclass[10pt]{amsart}
\usepackage{amsfonts}
\usepackage{ifthen}
\usepackage{amsthm}
\usepackage{amsmath,mathrsfs}
\usepackage{graphicx}
\usepackage{amscd,amssymb,amsthm}

\newcounter{minutes}
\divide\time by 60
\newcounter{hours}
\multiply\time by 60 \addtocounter{minutes}{-\time}

\newtheorem{theorem}{Theorem}

\newtheorem{corollary}{Corollary}

\newtheorem{problem}{Open problem}
\newtheorem{conjecture}{Conjecture}

\title{Zeros of Bessel function derivatives}

\author[\'A. Baricz]{\'Arp\'ad Baricz}
\address{Department of Economics, Babe\c{s}-Bolyai University, 400591 Cluj-Napoca, Romania}
\address{Institute of Applied Mathematics, \'Obuda University, 1034 Budapest, Hungary}
\email{bariczocsi@yahoo.com}

\author[C.G. Kokologiannaki]{Chrysi G. Kokologiannaki}
\address{Department of Mathematics, University of Patras, 26500 Patras, Greece}
\email{chrykok@math.upatras.gr}

\author[T.K. Pog\'any]{Tibor K. Pog\'any}
\address{Faculty of Maritime Studies, University of Rijeka, 51000 Rijeka, Croatia}
\address{Institute of Applied Mathematics, \'Obuda University, 1034 Budapest, Hungary}
\email{poganj@pfri.hr}



\begin{document}

\def\thefootnote{}
\footnotetext{ \texttt{File:~\jobname .tex,
          printed: \number\year-\number\month-\number\day,
          \thehours.\ifnum\theminutes<10{0}\fi\theminutes}
} \makeatletter\def\thefootnote{\@arabic\c@footnote}\makeatother

\maketitle

\begin{abstract}
We prove that for $\nu>n-1$ all zeros of the $n$th derivative of Bessel function of the first kind $J_{\nu}$ are real. Moreover, we show that the positive zeros of the $n$th and $(n+1)$th derivative of Bessel function of the first kind $J_{\nu}$ are interlacing when $\nu\geq n,$ and $n$ is a natural number or zero. Our methods include the Weierstrassian representation of the $n$th derivative, properties of the Laguerre-P\'olya class of entire functions, and the Laguerre inequalities. Some similar results for the zeros of the first and second derivative of the Struve function of the first kind $\mathbf{H}_{\nu}$ are also proved. The main results obtained in this paper generalize and complement some classical results on the zeros of Bessel and Struve functions of the first kind. Some conjectures and open problems related to Hurwitz theorem on the zeros of Bessel functions are also proposed, which may be of interest for further research.
\end{abstract}

\section{\bf Introduction and the Main Results}

The zeros of Bessel functions of the first kind are important in several problems of applied mathematics, mathematical physics and engineering sciences. Because of their importance there is an extensive literature on various properties of the zeros of Bessel functions of the first kind, and they were investigated by famous researchers like Bessel, Euler, Fourier, Gegenbauer, Hurwitz, Lommel, Rayleigh and Stokes. We refer to the survey paper \cite{kerimov} and to the references therein for detailed information on various properties of the zeros of Bessel functions of the first kind. In the last three decades the zeros of the $n$th derivative of Bessel functions of the first kind for $n\in\{1,2,3\}$ have been also studied by researchers like Elbert, Ifantis, Ismail, Kokologiannaki, Laforgia, Landau, Lorch, Mercer, Muldoon, Petropoulou, Siafarikas and Szeg\H o; for more details see the papers \cite{ismail,koko} and the references therein. However, the literature remains in silence about the zeros of the $n$th derivative of Bessel functions, when $n$ is a natural number greater than 4. Some interesting results and open problems were stated by Lorch and Muldoon in \cite{lorch} about $J_{\nu}^{(n)}$ when $\nu\in(n-3,n-2)$, however this is so far the only one study in which the zeros of higher order derivatives of $J_{\nu}$ were considered. In this paper our aim is to fill partially this gap and to present some results for the derivatives of Bessel functions by using the Laguerre-P\'olya class of entire functions and the so-called Laguerre inequalities. By using a similar technique as for the Bessel functions of the first kind our aim is also to present some new results for the zeros of the first and second derivative of the Struve functions of the first kind. Moreover, by using these results we find the explicit expressions for some Rayleigh sums for the zeros of Struve functions and their first and second derivatives. Also, lower bounds for the first positive zero of them are given. At the end of this section some conjectures and open problems are also proposed, which may be of interest for further research.

\subsection{Zeros of the $n$th derivative of Bessel functions} In view of the results on the zeros of the $n$th derivative of Bessel functions of the first kind, when $n\in\{0,1,2,3\},$ the statements of the following theorem are very natural and somehow expected, they provide the extensions of some classical results on the zeros of Bessel function of the first kind and its derivative of the first order. Throughout of this paper $n$ and $s$ are from $\mathbb{N}_0=\{0,1,2,\dots\}$.

\begin{theorem}\label{th1}
The following assertions are valid:
\begin{enumerate}
\item[\bf a.] If $\nu>n-1,$ then $x\mapsto J_{\nu}^{(n)}(x)$ has infinitely many zeros, which are all real and simple, except the origin.
\item[\bf b.] If $\nu\geq n,$ then the positive zeros of the $n$th and $(n+1)$th derivative of $J_{\nu}$ are interlacing.
\item[\bf c.] If $\nu>n-1,$ then all zeros of $x\mapsto (n-\nu)J_{\nu}^{(n)}(x)+xJ_{\nu}^{(n+1)}(x)$ are real and interlace with the zeros of $x\mapsto J_{\nu}^{(n)}(x).$
\end{enumerate}
\end{theorem}

It is important to mention that part {\bf b} in particular reduces to the chains of inequalities
$$j_{\nu,1}'''<j_{\nu,1}''<j_{\nu,2}'''<j_{\nu,2}''<j_{\nu,3}'''<j_{\nu,3}''<\dots,\ \ \ \nu\geq2$$
and
$$j_{\nu,1}''<j_{\nu,1}'<j_{\nu,2}''<j_{\nu,2}'<j_{\nu,3}''<j_{\nu,3}'<\dots,\ \ \ \nu\geq1,$$
where $j_{\nu,n}''$ and $j_{\nu,n}'''$ denote the $n$th positive zero of $J_{\nu}''$ and $J_{\nu}''',$ respectively. These inequalities complement the well-known ones \cite[p. 235]{nist}
$$j_{\nu,1}'<j_{\nu,1}<j_{\nu,2}'<j_{\nu,2}<j_{\nu,3}'<j_{\nu,3}<\dots,\ \ \ \nu\geq0.$$

We also note that part {\bf c} is actually a generalization of the well-known result that for $\nu>-1$ the zeros of the Bessel functions $J_{\nu}$ and $J_{\nu+1}$ are interlacing (see \cite[p. 235]{nist}). Namely, by choosing $n=0$ in part {\bf c} we get that the zeros of $J_{\nu}$ and $x\mapsto xJ_{\nu}'(x)-\nu J_{\nu}(x)=xJ_{\nu+1}(x)$ are interlacing.

Finally, it is worth to mention that, by using the main idea from \cite[p. 705]{dimitar}, an immediate consequence of part {\bf a} of Theorem \ref{th1} in terms of generalized hypergeometric polynomials reads as follows. The connection between these two results is a special class of real entire functions, called Laguerre-P\'olya class of entire functions. Recall that a real entire function $\psi$ belongs to the Laguerre-P\'{o}lya class $\mathcal{LP}$ if it can be represented in the form
$$\psi(x) = c x^{m} e^{-a x^{2} + b x} \prod_{n\geq1}
\left(1+\frac{x}{x_{n}}\right) e^{-\frac{x}{x_{n}}},$$
with $c,$ $b,$ $x_{n} \in \mathbb{R},$ $a \geq 0,$ $m\in
\mathbb{N}_0$ and $\sum_{n\geq1} x_{n}^{-2} < \infty.$ We note that the class $\mathcal{LP}$ consists of entire functions which are uniform limits
on the compact sets of the complex plane of polynomials with only real zeros. For more details on the class $\mathcal{LP}$ we refer to \cite[p. 703]{dimitar}
and to the references therein.

\begin{theorem}\label{th2}
If $\nu>n-1,$ then all the zeros of the Laguerre-type polynomial $${}_3F_3\left(-s,\frac{\nu+1}2, \frac\nu2 +1; \nu+1, \frac{\nu-n+1}2, \frac{\nu-n}2 +1; x\right)$$
are real and simple.
\end{theorem}

We note that the denomination Laguerre-type polynomial for the Jensen polynomial appearing in Theorem \ref{th2} is justified by the facts that the case $n=0$ reduces to the well-known generalized Laguerre polynomial ${}_1F_1(-s;\nu+1;x)$ (see \cite[p. 705]{dimitar}), the case $n=1$ one corresponds to the generalized hypergeometric polynomial ${}_2F_2 \left(-s,\frac\nu2 +1;\nu+1, \frac{\nu}2; x \right)$ which one transforms into the Koornwinder's generalized Laguerre polynomial \cite[p. 26]{Koekoek}, while the case $n=2$ is related to the generalized Laguerre polynomial ${}_3F_3\left(-s,\frac{\nu+1}2, \frac\nu2 +1; \nu+1, \frac{\nu-1}2, \frac{\nu}2; x\right),$ considered by \'Alvarez-Nodarse and Marcell\'an \cite{marcellan}.

\subsection{Zeros of the first and second derivatives of Struve functions} In 1970 Steinig \cite{steinig} studied the real zeros of the Struve functions, and proved that when $|\nu|<\frac{1}{2}$ the zeros of the Struve function $\mathbf{H}_{\nu}$ are all real and simple, the positive zeros of $\mathbf{H}_{\nu}$ interlace with the positive zeros of $J_{\nu},$ and lie in the intervals $(m\pi, (m+1)\pi),$ $m\in\mathbb{N}.$ Motivated by these results and in order to find the radius of convexity of some normalized Struve functions, recently Baricz and Ya\u{g}mur \cite{nihat} proved that the zeros of the function $\mathbf{H}_{\nu}'$ are all real and simple for $|\nu|<{1\over 2},$ and the positive zeros of the function $\mathbf{H}_{\nu}'$ interlace with the positive zeros of $\mathbf{H}_{\nu}.$ In this subsection we are going to prove some analogous results for the second derivative of $\mathbf{H}_{\nu}.$

\begin{theorem}\label{th3}
The following assertions are valid:
\begin{enumerate}
\item[\bf a.] If $\nu\in\left(0,\frac{1}{2}\right],$ then $x\mapsto \mathbf{H}_{\nu}''(x)$ has infinitely many zeros, which are all real and simple.
\item[\bf b.] If $\nu\in\left(0,\frac{1}{2}\right],$ then the positive zeros of the first and second derivative of $\mathbf{H}_{\nu}$ are interlacing.
\item[\bf c.] If $\nu\in\left[-\frac{1}{2},\frac{1}{2}\right],$ then all zeros of $x\mapsto -\nu\mathbf{H}_{\nu}'(x)+x\mathbf{H}_{\nu}''(x)$ are real and interlace with the zeros of $x\mapsto \mathbf{H}_{\nu}'(x).$ Moreover, all zeros of $x\mapsto -(\nu+1)\mathbf{H}_{\nu}(x)+x\mathbf{H}_{\nu}'(x)$ are real and interlace with the zeros of $x\mapsto \mathbf{H}_{\nu}(x).$
\end{enumerate}
\end{theorem}

The next result is the counterpart of Theorem \ref{th2} for Struve functions.

\begin{theorem}\label{th4}
If $\nu\in\left[-\frac{1}{2},\frac{1}{2}\right]$ and $n\in\{0,1\},$ or $\nu\in\left(0,\frac{1}{2}\right]$ and $n=2,$ then all the zeros of the hypergeometric
polynomial
   \[ {}_4F_4 \left(-s,1, \frac\nu2+1,\frac\nu2+\frac32;\frac32,\nu + \frac32,\frac{\nu-n}2+1,\frac{\nu-n}2+\frac32; x\right)\]
are real and simple.
\end{theorem}

\subsection{Rayleigh sums of the zeros of the first and second derivative of Struve functions} The Struve function $\mathbf{H}_{\nu}$ and its derivative of the first and second order can be represented by infinite series as well as using Hadamard's factorization as by infinite products. So, equating these representations, as it was done in \cite{KL} for the zeros of $J_{\nu}$ and $J_{\nu}',$ it is possible to obtain the Rayleigh sums for the zeros of the Struve function $\mathbf{H}_{\nu}$ and of its first and second derivatives. For more information on the Rayleigh sums of the zeros of Bessel functions of the first kind we refer to \cite[p. 240]{nist}, \cite[p. 502]{watson} and to the references therein.

\begin{theorem}\label{th5}
If $\nu\in\left[-\frac{1}{2},\frac{1}{2}\right],$ then the first two Rayleigh sums for the $n$th positive zeros $h_{\nu,n},$ $h_{\nu,n}'$ and $h_{\nu,n}''$ of the Struve functions $\mathbf{H}_{\nu},$ $\mathbf{H}_{\nu}'$ and $\mathbf{H}_{\nu}''$ are
$$\sum_{n\geq 1}{1\over {h_{\nu,n}^{2}}}={1\over {3(2\nu+3)}},\ \ \sum_{n\geq 1}{1\over {h_{\nu,n}^{4}}}={{7-2\nu}\over {45(2\nu+3)^{2}(2\nu+5)}},$$
$$\sum_{n\geq 1}{1\over {(h_{\nu,n}^{\prime})^{2}}}={{\nu+3}\over {3(\nu+1)(2\nu+3)}},\ \ \sum_{n\geq 1}{1\over {(h_{\nu,n}^{\prime})^{4}}}={{-2\nu^{3}-5\nu^{2}+72\nu+135}\over {45(\nu+1)^{2}(2\nu+3)^2(2\nu+5)}},
$$
$$
\sum_{n\geq 1}{1\over {(h_{\nu,n}'')^{2}}}={{(\nu+2)(\nu+3)}\over {3\nu(\nu+1)(2\nu+3)}},\
\sum_{n\geq 1}{1\over {(h_{\nu,n}'')^{4}}}={{-2\nu^{5}-13\nu^{4}+92\nu^{3}+763\nu^{2}+1500\nu+900}\over {45\nu^{2}(\nu+1)^{2}(2\nu+3)^{2}(2\nu+5)}},
$$
provided that $\nu>0$ in the last two relations.
\end{theorem}

An immediate consequence of the above theorem is the following result. We note that if we keep on doing the same procedure as in the proof of Theorem \ref{th5}, then we can derive the sums $\sum_{n\geq 1}{{h_{\nu,n}^{-2k}}}$, $k\in\{3,4,\dots\}$ and by using these Rayleigh sums it is possible to derive sharper lower bounds for $h_{\nu,1},$ $h_{\nu,1}'$ and $h_{\nu,1}''.$

\begin{corollary}\label{cor}
For $|\nu|<\frac{1}{2}$ we have the following inequalities:
$$h_{\nu,1}>\sqrt{3(2\nu+3)},\ \ h_{\nu,1}^{2}>3(2\nu+3)\sqrt{{{5(2\nu+5)}\over{7-2\nu}}},$$
$$h_{\nu,1}^{\prime}>\sqrt{{{3(\nu+1)(2\nu+3)}\over{\nu+3}}},\ \ \
(h_{\nu,1}^{\prime})^{2}>3(\nu+1)(2\nu+3)\sqrt{\frac{5(2\nu+5)}{-2\nu^3-5\nu^2+72\nu+135}},$$
$$h_{\nu,1}^{\prime\prime}>\sqrt{{{3\nu(\nu+1)(2\nu+3)}\over{(\nu+2)(\nu+3)}}}, \ \
(h_{\nu,1}^{\prime\prime})^{2}>{{3\nu(\nu+1)(2\nu+3)}}\sqrt{{{5(2\nu+5)}\over{-2\nu^{5}-13\nu^{4}+92\nu^{3}+763\nu^{2}+1500\nu+900}}},$$
provided that $\nu>0$ in the last two inequalities.
\end{corollary}

Since $h_{\nu,1}^{\prime}<h_{\nu,1}$, we can see that the third and fourth inequalities in the above corollary are also bounds for $h_{\nu,1},$ but sharper than the ones given in the first and second inequalities. Moreover, for $\nu\in\left(0,\frac{1}{2}\right]$ we have that $h_{\nu,1}''<h_{\nu,1}'$ and thus the fifth and sixth inequalities in the above corollary are also bounds for $h_{\nu,1}'$ and $h_{\nu,1}.$

\subsection{Conjectures and open problems concerning the zeros of derivatives of Bessel and Struve functions} Part {\bf a} of Theorem \ref{th1} is an extension of the celebrated result of von Lommel (see \cite[p. 482]{watson}), which states that if $\nu>-1$ then all zeros of $J_{\nu}$ are real. Now, recall that by means of Hurwitz's theorem (see \cite{hu}, \cite[p. 483]{watson}) we know that if $\nu>-1,$ then all zeros of $J_{\nu}$ are real; when $-2s-2<\nu<-2s-1,$ $s\in\mathbb{N}_0,$ then $J_{\nu}$ has $4s+2$ complex zeros, of which two are purely imaginary; and when $-2s-1<\nu<-2s,$ $s\in\mathbb{N},$ then the Bessel function $J_{\nu}$ has $4s$ complex zeros, of which none are purely imaginary. See also \cite{hu,kerimov,kim} for more details. Thus, it is an interesting conjecture to provide a complete picture about the behavior of the zeros of the derivatives of Bessel functions, that is, to present a generalization of Hurwitz theorem \cite{hu}.

\begin{conjecture}
Let $n\in\mathbb{N}_0.$
\begin{enumerate}
\item[\bf a.] Is it true that if $s$ is a nonnegative integer and $n-2s-2<\nu<n-2s-1,$ then the function $J_{\nu}^{(n)}$ has $4s+2$ complex zeros, of which two are purely imaginary?
\item[\bf b.] Is it true that if $s$ is a positive integer and $n-2s-1<\nu<n-2s,$ then the function $J_{\nu}^{(n)}$ has $4s$ complex zeros, of which none are purely imaginary?
\end{enumerate}
\end{conjecture}

Note that Hurwitz proof of his famous theorem on the distribution of the zeros of Bessel functions of the first kind was based on the three term recurrence relation of Bessel functions, on Lommel polynomials and on the relation between Bessel functions and Lommel polynomials, see \cite{hu,kerimov} and \cite[p. 483]{watson} for more details. Sixteen years ago Ki and Kim \cite{kim} presented a nice alternative proof of Hurwitz theorem by using a Fourier critical point approach. To define the notion of the Fourier critical point let $f$ be a real entire function defined in an open interval $(a,b)\subset\mathbb{R}.$ Let $l\in\mathbb{N}$ and suppose that $c\in(a,b)$ is a zero of $f^{(l)}(x)$ of multiplicity $m\in\mathbb{N},$ that is, $f^{(l)}(c)=\dots=f^{(l+m-1)}(c)=0$ and $f^{(l+m)}(c)\neq0.$ Now, let $k=0$ if $f^{(l-1)}(c)=0,$ otherwise let
$$k=\left\{\begin{array}{ll}{m}/{2},& \mbox{if}\ m\ \mbox{is even},\\
{(m+1)}/{2},& \mbox{if}\ m\ \mbox{is odd and}\ f^{(l-1)}(c)f^{(l+m)}(c)>0,\\
{(m-1)}/{2},& \mbox{if}\ m\ \mbox{is odd and}\ f^{(l-1)}(c)f^{(l+m)}(c)<0. \end{array}\right.$$
Then we say that $f^{(l)}(x)$ has $k$ critical zeros and $m-k$ noncritical zeros at $x = c.$ A point in $(a,b)$ is said to be a Fourier critical point of $f$ if some derivative of $f$ has a critical zero at the point. For more details on Fourier critical points we refer to \cite{kim}.

Motivated by the proof of Hurwitz theorem given in \cite{kim}, we consider the auxiliary function
$$f_{\nu,n}(x)=\sum_{m\geq0}\frac{\Gamma(\nu+2m+1)}{\Gamma(\nu+2m-n+1)\Gamma(\nu+m+1)}\frac{x^m}{m!}.$$
This function is a real entire function of growth order $\frac{1}{2}$ and consequently of genus $0$ and according to \cite[Theorem 4.1]{kim} has just as many Fourier critical points as couples of nonreal zeros. Now, since $2^nx^{{n}/{2}}J_{\nu}^{(n)}(2\sqrt{x})=x^{{\nu}/{2}}f_{\nu,n}(-x),$
part {\bf a} of Theorem \ref{th1} implies that when $\nu>n-1$ the function $f_{\nu,n}$ has no Fourier critical points. The following conjecture is motivated by \cite[p. 68]{kim} and in the case that it would be true it would imply that the answers for the questions stated in the first conjecture are affirmative.

\begin{conjecture}
Let $n\in\mathbb{N}_0.$
\begin{enumerate}
\item[\bf a.] Is it true that if $s$ is a nonnegative integer and $n-2s-2<\nu<n-2s-1,$ then the function $f_{\nu,n}$ has exactly $s$ Fourier critical points and one positive real zero?
\item[\bf b.] Is it true that if $s$ is a positive integer and $n-2s-1<\nu<n-2s,$ then the function $f_{\nu,n}$ has exactly $s$ Fourier critical points and no positive real zeros?
\end{enumerate}
\end{conjecture}

Owing to Steinig \cite[p. 367]{steinig} we know that when $\nu\in\left[-\frac{1}{2},\frac{1}{2}\right]$ all zeros of the Struve function $\mathbf{H}_{\nu}$ are real. We also know that when $\nu>\frac{1}{2}$ we have that $\mathbf{H}_{\nu}(x)>0$ for $x>0,$ and thus there are no positive real zeros in this case. However, we do not know what happens with the zeros when $\nu<-\frac{1}{2}.$ By using the connection between Bessel and Struve functions, that is, the relation $\mathbf{H}_{-m-\frac{1}{2}}(x)=(-1)^mJ_{m+\frac{1}{2}}(x),$ where $m\in\mathbb{N}_0,$ it is clear that all zeros of the Struve function $\mathbf{H}_{\nu}$ are real when $\nu=-m-\frac{1}{2}$ and $m\in\mathbb{N}_0.$ Moreover, the same is true for $\mathbf{H}_{\nu}',$ and in general we can state that if $\nu=-m-\frac{1}{2},$ $m+\frac{1}{2}>n-1$ and $m,n\in\mathbb{N}_0,$ then all zeros of the function $\mathbf{H}_{\nu}^{(n)}$ are real. The next two open problems are motivated by this result.

\begin{problem}
Find the number of complex zeros of the function $\mathbf{H}_{\nu}^{(n)}$ when $\nu\in\mathbb{R}$ and $n\in\mathbb{N}_0.$
\end{problem}

\begin{problem}
Is it true that if $\nu\in\left[-\frac{1}{2},\frac{1}{2}\right],$ and $n\in\mathbb{N}_0,$ then all zeros of $\mathbf{H}_{\nu}^{(n)}$ are real? Find the maximal range for $\nu$ for which all zeros of $\mathbf{H}_{\nu}^{(n)}$ are real.
\end{problem}

Recall that Steinig showed in \cite[p. 367]{steinig} that the zeros of the functions $J_{\nu}$ and $\mathbf{H}_{\nu}$ are interlacing when $|\nu|<\frac{1}{2}.$ The last open problem of this paper is motivated by this result.

\begin{problem}
Is it true that if $\nu\in\left[-\frac{1}{2},\frac{1}{2}\right],$ and $n\in\mathbb{N}_0,$ then the zeros of $\mathbf{H}_{\nu}^{(n)}$ and $J_{\nu}^{(n)}$ interlace? Find the maximal range for $\nu$ for which all zeros of $\mathbf{H}_{\nu}^{(n)}$ and $J_{\nu}^{(n)}$ interlace.
\end{problem}

Recently, Baricz and Sz\'asz \cite{bs}, and Baricz et al. \cite{bcd} found necessary and sufficient conditions on the parameter $\nu$ such that for $n\in\{0,1,2,3\}$ the function $z\mapsto 2^{\nu}\Gamma(\nu-n+1)z^{\frac{n+2-\nu}{2}}J_{\nu}^{(n)}(\sqrt{z})$ is starlike (maps the open unit disk of the complex plane into
a starlike domain) and all of its derivatives are close-to-convex (and hence univalent). In the proofs the key tool was that for fixed $m\in\mathbb{N}$ and $n\in\mathbb{N}_0$ the $m$th positive zeros of $J_{\nu}^{(n)},$ denoted by $j_{\nu,m}^{(n)},$ are increasing with $\nu$ on $(n-1,\infty),$ where $n\in\{0,1,2,3\}$ (see \cite{koko,lorch2,mercer,nist,wong} for more details). Now, since the zeros of $J_{\nu}^{(n)}$ for $\nu>n-1$ are real, an affirmative answer to the following conjecture would enable to generalize the above exposed results of \cite{bs} and \cite{bcd}.

\begin{conjecture}
Is it true that $\nu\mapsto j_{\nu,m}^{(n)}$ is increasing on $(n-1,\infty)$ for $n\in\mathbb{N}_0$ and $m\in\mathbb{N}$ fixed?
\end{conjecture}

\section{\bf Proof of the main results}

\begin{proof}[\bf Proof of Theorem \ref{th1}]
{\bf a.} We shall use mathematical induction to prove the reality of the zeros. When $n=\{0,1,2,3\},$ then we know that for $\nu>n-1$ the zeros of $J_{\nu}^{(n)}(x)$ are all real, see for example \cite{ismail,kerimov,koko} and the references therein for more details. We suppose that for $n\in\{4,5,\dots\}$ fixed and $\nu>n-1$ we have that $J_{\nu}^{(n)}(x)$ has only real zeros and we show that when $\nu>n$ then $J_{\nu}^{(n+1)}(x)$ has also only real zeros. We denote the $m$th positive zero of $J_{\nu}^{(n)}(x)$ by $j_{\nu,m}^{(n)},$ where $m\in\mathbb{N}$ and $n\in\mathbb{N}_0.$ Due to Skelton \cite[p. 340]{skelton} we know that the following Weierstrassian decomposition holds
\begin{equation}\label{prod}J_{\nu}^{(n)}(x)=\frac{x^{\nu-n}}{2^{\nu}\Gamma(\nu+1-n)}\prod_{m\geq 1}\left(1-\frac{x^2}{\left(j_{\nu,m}^{(n)}\right)^2}\right),\end{equation}
and this infinite product is uniformly convergent on compact subsets of the complex plane. Note that it was not stated in \cite{skelton} for which $\nu$ is valid the above infinite product. Namely, \cite[Theorem 2.1]{skelton} is enounced for $\nu\geq n,$ and the above infinite product \eqref{prod} appears in the proof of \cite[Theorem 2.1]{skelton}. However, \eqref{prod} holds true for all $\nu>n-1,$ where $n$ is a natural number or zero. To see this we note that since for fixed $n$ and $\nu$
$$\lim_{m\to\infty}\frac{m\log m}{\log\frac{2^{2m}}{\Gamma(\nu+1-n)}+\log\Gamma(m+1)+\log\Gamma(\nu+m+1)+\log\Gamma(\nu+2m-n+1)-\log\Gamma(\nu+2m+1)}=\frac{1}{2},$$
the real entire function
$$x\mapsto \mathbb{J}_{\nu,n}(x)=2^{\nu}\Gamma(\nu+1-n)x^{n-\nu}J_{\nu}^{(n)}(x)=
\sum_{m\geq0}\frac{(-1)^m\Gamma(\nu+2m+1)\Gamma(\nu+1-n)x^{2m}}{m!2^{2m}\Gamma(\nu+2m-n+1)\Gamma(\nu+m+1)}$$
is of order $\frac{1}{2}$ and in view of the Hadamard theorem \cite[p. 26]{lev} it follows that \eqref{prod} is indeed true for $\nu>n-1$. Here we used the limit $\frac{\log\Gamma(am+b)}{m\log m}\to a,$ as $m\to\infty,$ where $a,b>0.$ To see this just observe that
$$\lim_{x\to\infty}\frac{\log\Gamma(ax+b)}{x\log x}=a\lim_{x\to\infty}\frac{\psi(ax+b)}{1+\log x}=a\lim_{x\to\infty}\frac{\log(ax+b)-\frac{1}{2(ax+b)}+ \mathcal{O}(x^{-2})}{\log x}=a,$$
where $\psi(x)=\Gamma'(x)/\Gamma(x)$ is the logarithmic derivative of the Euler gamma function. On the other hand, it is known that (see \cite[Theorem 2.9.2]{boas}) if $f$ is an entire function and its growth order is finite and is not equal to a positive integer, then $f$ has infinitely many zeros, or $f$ is a polynomial. Thus, using the fact that the growth order of the real entire function $x\mapsto x^{n-\nu}J_{\nu}^{(n)}(x)$ is $\frac{1}{2},$ and by using the above mentioned result we obtain that $J_{\nu}^{(n)}(x)$ has infinitely many zeros when $\nu>n-1.$ Moreover, by using \eqref{prod} we obtain
$$J_{\nu}^{(n+1)}(x) = \frac{x^{\nu-n} \prod\limits_{m\geq 1}\left(1-\frac{x^2}{\left(j_{\nu,m}^{(n)}\right)^2}\right)}{2^\nu \Gamma(\nu+1-n)}\left(\frac{\nu-n}{x}- \sum_{m \geq 1} \frac{2x}{\left(j_{\nu,m}^{(n)}\right)^2 - x^2}\right),$$
which in turn implies that
\begin{equation}\label{nml}
\frac{J_{\nu}^{(n+1)}(x)}{J_{\nu}^{(n)}(x)}=\frac{\nu-n}{x}- \sum_{m \geq 1} \frac{2x}{\left(j_{\nu,m}^{(n)}\right)^2 - x^2},
\end{equation}
where $\nu>n,$ $x$ is real (or complex) such that $x\neq0,$ $x\neq j_{\nu,m}^{(n)},$ $m\in\mathbb{N}$ and $n\in\mathbb{N}_0.$ Now, we are going to conclude by induction on $n$ that for $\nu>n$ all the zeros of $J_\nu^{(n+1)}(x)$ are real, provided that for $\nu>n-1$ all the zeros of $J_\nu^{(n)}(x)$ are real. For this first we show that for $\nu>n$ the zeros of $J_\nu^{(n+1)}(x)$ cannot be purely imaginary. Indeed, letting $J_\nu^{(n+1)}({\rm i}y) = 0,$ where $y \in \mathbb{R},$ $y\neq0,$ from \eqref{nml} we deduce
   \[ 0 = -{\rm i}\left(\nu-n +  \sum_{m \geq 1} \frac{2y^2}{\left(j_{\nu,m}^{(n)}\right)^2 +y^2}\right),\]
which is a contradiction since $\nu>n$ and the zeros $j_{\nu,m}^{(n)}$ are real. Finally, we show that for $\nu>n$ the zeros of $J_\nu^{(n+1)}(x)$ cannot be complex. Taking $z = x+{\rm i}y,$ where  $xy\neq 0,$ a complex zero of $J_\nu^{(n+1)}(z)$ and denoting by $\omega = \left(j_{\nu,m}^{(n)}\right)^2-x^2+y^2$, in view of \eqref{nml} we have
   \[ z\frac{J_{\nu}^{(n+1)}(z)}{J_{\nu}^{(n)}(z)} = \nu - n - 2\sum_{m \geq 1} \frac{(x^2-y^2)\omega-4x^2y^2}{\omega^2 + 4x^2y^2}
			      - 4{\rm i}xy \sum_{m \geq 1} \frac{\omega+x^2-y^2}{\omega^2 + 4x^2y^2} = 0,\]
that is
$$\sum_{m \geq 1} \frac{\omega+x^2-y^2}{\omega^2 + 4x^2y^2} = \sum_{m \geq 1} \frac{\left(j_{\nu,m}^{(n)}\right)^2}{\omega^2 + 4x^2y^2} = 0,$$
which is a contradiction. Thus, indeed for $\nu>n$ the zeros of $J_{\nu}^{(n+1)}(x)$ are all real.

Now, we prove that these zeros are all simple, except the origin. If we suppose that $\rho\neq0$ is a double zero of $J_{\nu}^{(n)}(z)$ it follows that the derivative of the quotient $J_{\nu}^{(n)}(z)/J_{\nu}^{(n-1)}(z)$ also vanishes at $\rho,$ which is contradiction since in view of \eqref{nml} for $\nu>n-1$ we have that
$$\frac{d}{dz}\left(\frac{J_{\nu}^{(n)}(z)}{J_{\nu}^{(n-1)}(z)}\right)=-\frac{\nu-n+1}{z^2}- 2\sum_{m \geq 1} \frac{\left(j_{\nu,m}^{(n-1)}\right)^2 + z^2}{\left(\left(j_{\nu,m}^{(n-1)}\right)^2 - z^2\right)^2}\neq0.$$

{\bf b.} Since the zeros of $J_{\nu}^{(n)}(x)$ are all real, it follows that the function $\mathbb{J}_{\nu,n}$ is in the Laguerre-P\'olya class of real entire functions since the exponential factors in the infinite product are canceled because of the symmetry of the zeros $\pm j_{\nu,m}^{(n)},$ $m\in\mathbb{N},$ with respect to the origin. Now, since $\mathbb{J}_{\nu,n}\in\mathcal{LP},$ it follows that it satisfies the Laguerre inequality (see \cite[p. 67]{skov})
$$\left(\mathbb{J}_{\nu,n}^{(k)}(x)\right)^2-\mathbb{J}_{\nu,n}^{(k-1)}(x)\mathbb{J}_{\nu,n}^{(k+1)}(x)>0,$$
where $n\in\mathbb{N}_0,$ $k\in\mathbb{N},$ $\nu>n-1$ and $x\in\mathbb{R}.$ Choosing $k=1$ in the above inequality we get
$$\left(x{J}_{\nu}^{(n+1)}(x)\right)^2-x^2{J}_{\nu}^{(n+2)}(x){J}_{\nu}^{(n)}(x)+(n-\nu)\left({J}_{\nu}^{(n)}(x)\right)^2>0,$$
which implies that
$$\left({J}_{\nu}^{(n+1)}(x)\right)^2-{J}_{\nu}^{(n+2)}(x){J}_{\nu}^{(n)}(x)>(\nu-n)\left({J}_{\nu}^{(n)}(x)\right)^2/x^2>0,$$
where $\nu>n\geq0$ and $x\neq0.$ Consequently, the function $x\mapsto {J}_{\nu}^{(n+1)}(x)/{J}_{\nu}^{(n)}(x)$ is strictly decreasing on each interval $\left(j_{\nu,m-1}^{(n)},j_{\nu,m}^{(n)}\right),$ $m\in\mathbb{N}$ (note that $j_{\nu,0}^{(n)}=0$). On the other hand, for fixed $m\in\mathbb{N}$ the function $x\mapsto {J}_{\nu}^{(n+1)}(x)/{J}_{\nu}^{(n)}(x)$ takes the limit $\infty$ when $x\searrow j_{\nu,m-1}^{(n)},$ and the limit $-\infty$ when $x\nearrow j_{\nu,m}^{(n)}.$ Summarizing, for arbitrary $m\in\mathbb{N}$ the restriction of the function $x\mapsto {J}_{\nu}^{(n+1)}(x)/{J}_{\nu}^{(n)}(x)$ in each interval $\left(j_{\nu,m-1}^{(n)},j_{\nu,m}^{(n)}\right)$ intersects the horizontal line only once, and the abscissa of this intersection point is exactly $j_{\nu,m}^{(n+1)}.$ With this we proved that when $\nu>n$ the positive zeros of $J_{\nu}^{(n+1)}(x)$ and $J_{\nu}^{(n)}(x)$ are interlacing.

It is worth to mention that the monotonicity of $x\mapsto {J}_{\nu}^{(n+1)}(x)/{J}_{\nu}^{(n)}(x)$ can be also verified by using the Mittag-Leffler expansion \eqref{nml}. Namely, we have
$$\frac{d}{dx}\left(\frac{J_{\nu}^{(n+1)}(x)}{J_{\nu}^{(n)}(x)}\right)=-\frac{\nu-n}{x^2}- 2\sum_{m \geq 1} \frac{\left(j_{\nu,m}^{(n)}\right)^2 + x^2}{\left(\left(j_{\nu,m}^{(n)}\right)^2 - x^2\right)^2}<0$$
for all $\nu>n,$ $n\in\mathbb{N}_0$ and $x\neq j_{\nu,m}^{(n)},$ $m\in\mathbb{N}_0.$

{\bf c.} Recall that Laguerre's theorem on separation of zeros \cite[p. 23]{boas} states that, if $z\mapsto f(z)$ is an entire function, not a constant,
which is real for real $z$ and has only real zeros, and is of genus $0$ or $1,$ then the zeros of $f'$ are also real and are separated by the zeros of $f.$ Now, according to the proof of part {\bf a} we get that $\mathbb{J}_{\nu,n}$ is a real entire function of genus zero. Thus, by using part {\bf a} of this theorem and Laguerre's separation theorem it results that the zeros of $x\mapsto (n-\nu)J_{\nu}^{(n)}(x)+xJ_{\nu}^{(n+1)}(x)$ are real when $\nu>n-1$ and are interlacing with the zeros of $x\mapsto J_{\nu}^{(n)}(x).$
\end{proof}

\begin{proof}[\bf Proof of Theorem \ref{th2}]
By means of \eqref{prod} and Theorem \ref{th1} the function
$$x\mapsto \mathbb{J}_{\nu,n}(2\sqrt{x})=2^n\Gamma(\nu+1-n)x^{\frac{n-\nu}{2}}J_{\nu}^{(n)}(2\sqrt{x})=
\sum_{m\geq0}\frac{(-1)^m\Gamma(\nu+2m+1)\Gamma(\nu+1-n)x^{m}}{m!\Gamma(\nu+2m-n+1)\Gamma(\nu+m+1)}$$ belongs also
to the Laguerre-P\'olya class $\mathcal{LP}.$ Consequently by using the well-known theorem of Jensen (see \cite{jensen} or \cite[Theorem A]{dimitar}) it follows that
the Jensen polynomial of $x\mapsto \mathbb{J}_{\nu,n}(2\sqrt{x})$ has only real zeros. Now, the Jensen polynomial in the question is
$$\sum_{m=0}^s (-1)^m \binom sm \frac{\Gamma(\nu+2m+1)\Gamma(\nu+1-n)}{\Gamma(\nu+2m-n+1)\Gamma(\nu+m+1)} x^m,$$
which after some transformations and in view of the Legendre duplication formula
\begin{equation}\label{eqL}\Gamma(2x) \sqrt{\pi} = 2^{2x-1} \Gamma(x)
\Gamma\left(x+\tfrac12\right)\end{equation} can be rewritten as
$${}_3F_3\left(-s,\frac{\nu+1}2, \frac\nu2 +1; \nu+1, \frac{\nu-n+1}2, \frac{\nu-n}2 +1; x\right).$$
Moreover, according to Csordas and Williamson \cite{csordas} the zeros of the Jensen polynomials are simple, and this completes the proof of the theorem.
\end{proof}

\begin{proof}[\bf Proof of Theorem \ref{th3}]
{\bf a.} We know that for $\nu\in\left[-\frac{1}{2},\frac{1}{2}\right]$ the zeros of $\mathbf{H}_{\nu}'(x)$ are all real and simple, see \cite{nihat}. For convenience we denote the $m$th positive zero of $\mathbf{H}_{\nu}^{(n)}(x)$ by $h_{\nu,m}^{(n)},$ where $m\in\mathbb{N}$ and $n\in\mathbb{N}_0.$ Since
$$\lim_{m\to\infty}\frac{m\log m}{\log2^{2m}+\log\Gamma\left(m+\frac{3}{2}\right)+\log\Gamma\left(m+\nu+\frac{3}{2}\right)-\log\left|(2m+\nu+1)(2m+\nu){\dots}(2m+\nu-n+2)\right|}=\frac{1}{2},$$
the real entire function
\begin{equation}\label{sumS}x\mapsto \mathbb{H}_{\nu,n}(x)=2^{\nu+1}x^{n-\nu-1}\mathbf{H}_{\nu}^{(n)}(x)=
\sum_{m\geq0}\frac{(-1)^m\Gamma(\nu+2m+1)x^{2m}}{2^{2m}\Gamma\left(\nu+\frac{3}{2}\right)\Gamma(\nu+2m-n+2)\Gamma\left(\nu+m+\frac{3}{2}\right)}\end{equation}
is of order $\frac{1}{2}$ and in view of the Hadamard theorem \cite[p. 26]{lev} it follows that the following Weierstrassian decomposition is valid for appropiate values of $\nu$ (for example $\nu\in\left[-\frac{1}{2},\frac{1}{2}\right],$ $\nu\neq0$) and $n\in\mathbb{N}$
\begin{equation}\label{prodS}\mathbf{H}_{\nu}^{(n)}(x)=\frac{(\nu+1)\nu\cdot{\dots}\cdot(\nu-n+2)x^{\nu+1-n}}{\sqrt{\pi}2^{\nu}\Gamma\left(\nu+\frac{3}{2}\right)}
\prod_{m\geq 1}\left(1-\frac{x^2}{(h_{\nu,m}^{(n)})^2}\right).\end{equation}
Thus, $\mathbf{H}_{\nu}^{(n)}(x)$ has infinitely many zeros and by using the fact that the zeros of $\mathbf{H}_{\nu}^{(n)}(x)$ are all real (and simple) when $n=1$, it follows that for $n=1$ the function $\mathbb{H}_{\nu,n}$ is in the Laguerre-P\'olya class of real entire functions since the exponential factors in the infinite product are canceled because of the symmetry of the zeros $\pm h_{\nu,m}^{(n)},$ $m\in\mathbb{N},$ with respect to the
origin. Now, since for $n=1$ we have $\mathbb{H}_{\nu,n}\in\mathcal{LP},$ it follows that it satisfies the Laguerre inequality (see \cite[p. 67]{skov})
$$\left(\mathbb{H}_{\nu,n}^{(k)}(x)\right)^2-\mathbb{H}_{\nu,n}^{(k-1)}(x)\mathbb{H}_{\nu,n}^{(k+1)}(x)>0,$$
where $n=1,$ $k\in\mathbb{N},$ $\nu\in\left[-\frac{1}{2},\frac{1}{2}\right]$ and $x\in\mathbb{R}.$ Choosing $k=1$ in the above inequality we get
$$\left(x\mathbf{H}_{\nu}^{(n+1)}(x)\right)^2-x^2\mathbf{H}_{\nu}^{(n+2)}(x)\mathbf{H}_{\nu}^{(n)}(x)+(n-\nu-1)\left(\mathbf{H}_{\nu}^{(n)}(x)\right)^2>0,$$
which implies that
$$\left(\mathbf{H}_{\nu}^{(n+1)}(x)\right)^2-\mathbf{H}_{\nu}^{(n+2)}(x)\mathbf{H}_{\nu}^{(n)}(x)>(\nu+1-n)\left(\mathbf{H}_{\nu}^{(n)}(x)\right)^2/x^2>0,$$
where $\nu\in\left(0,\frac{1}{2}\right]$ and $x\neq0.$ Consequently, for $n=1$ the function $x\mapsto \mathbf{H}_{\nu}^{(n+1)}(x)/\mathbf{H}_{\nu}^{(n)}(x)$ is strictly decreasing on each interval $\left(h_{\nu,m-1}^{(n)},h_{\nu,m}^{(n)}\right),$ $m\in\mathbb{N}.$ Here we used that $h_{\nu,0}^{(n)}=0.$ Since the zeros of $\mathbf{H}_{\nu}'$ are simple, the function $\mathbf{H}_{\nu}''$ does not vanish in $h_{\nu,m}'.$ On the other hand, when $n=1$ and for fixed $m\in\mathbb{N}$ the function $x\mapsto \mathbf{H}_{\nu}^{(n+1)}(x)/\mathbf{H}_{\nu}^{(n)}(x)$ takes the limit $\infty$ when $x\searrow h_{\nu,m-1}^{(n)},$ and the limit $-\infty$ when $x\nearrow h_{\nu,m}^{(n)}.$ Summarizing, when $n=1$ for arbitrary $m\in\mathbb{N}$ the restriction of the function $x\mapsto \mathbf{H}_{\nu}^{(n+1)}(x)/\mathbf{H}_{\nu}^{(n)}(x)$ in each interval $\left(h_{\nu,m-1}^{(n)},h_{\nu,m}^{(n)}\right)$ intersects the horizontal line only once, and the abscissa of this intersection point is exactly $h_{\nu,m}^{(n+1)}.$ Moreover, it is clear that these zeros are simple because of the above monotonicity and limit properties. Concerning the distribution of the zeros, a similar procedure shows that on the semi-axis $(-\infty,0)$ we have a similar situation as on the semi-axis $(0,\infty)$ and thus we proved that when $\nu\in\left(0,\frac{1}{2}\right]$ the zeros of $\mathbf{H}_{\nu}''(x)$ are all real and simple. With this the proof of part {\bf a} is complete, moreover, we also proved the statement of part {\bf b}.

{\bf c.} According to the proof of the previous part we get that $\mathbb{H}_{\nu,n}$ is a real entire function of genus zero. Thus, by using the fact that the zeros of $\mathbf{H}_{\nu}^{(n)}$ are all real when $n\in\{0,1\}$ and $\nu\in\left[-\frac{1}{2},\frac{1}{2}\right],$ in view of the Laguerre separation theorem it results that the zeros of $x\mapsto (n-\nu-1)\mathbf{H}_{\nu}^{(n)}(x)+x\mathbf{H}_{\nu}^{(n+1)}(x)$ are real when $\nu\in\left[-\frac{1}{2},\frac{1}{2}\right]$ and $n\in\{0,1\},$ and are interlacing with the zeros of $x\mapsto \mathbf{H}_{\nu}^{(n)}(x).$
\end{proof}

\begin{proof}[\bf Proof of Theorem \ref{th4}]
The $n$th derivative with respect to the argument of the Struve function of the order $\nu$ is the power series
   \begin{align*}
	    \mathbf H_\nu^{(n)}(x) &= \frac{x^{\nu-n+1}}{2^{\nu+1}} \sum_{m \geq 0} \frac{(-1)^m \Gamma(2m+\nu+2)\, x^{2m}}
			                          {2^{2m} \Gamma(m+\frac32)\Gamma(m+\nu+\frac32)\Gamma(2m+\nu-n+2)} \nonumber \\
														 &= \frac{x^{\nu-n+1}}{2^{\nu-n+1}} \sum_{m \geq 0}
														    \frac{(-1)^m \Gamma(m+\frac\nu2+1) \Gamma(m+\frac\nu2+\frac32)\, x^{2m}}
			                          {2^{2m} \Gamma(m+\frac32)\Gamma(m+\nu+\frac32)\Gamma(m+\frac{\nu-n}2+1)
																\Gamma(m+\frac{\nu-n}2+\frac32)}\,,
	 \end{align*}
where the duplication formula \eqref{eqL} is employed in both numerator and denominator. This means that we have accordingly
$$\mathbf H_\nu^{(n)}(2 \sqrt{x}) = \frac{2^{1-n}\,\Gamma(\nu+2)\, x^{\frac{\nu-n+1}2}}{\sqrt{\pi}\,
			                      \Gamma(\nu+\frac32) \Gamma(\nu-n+2)} \sum_{m \geq 0}
														\frac{(-1)^m (1)_m (\frac\nu2+1)_m (\frac\nu2+\frac32)_m}
			                      {(\frac32)_m (\nu+\frac32)_m (\frac{\nu-n}2+1)_m (\frac{\nu-n}2+\frac32)_m}\, \frac{x^m}{m!}\,$$
where $\nu-n+2\neq\{0,-1,\dots\},$ and the auxiliary function
\begin{align*} x \mapsto \mathcal H_{\nu, n}(2\sqrt{x}) &= \frac{\sqrt{\pi}\,\Gamma(\nu+\frac32)}{2^{1-n} \,\Gamma(\nu+2)}
			                                           \,\Gamma(\nu-n+2) \,x^{\frac{n-\nu-1}2}\, \mathbf H_\nu^{(n)}(2 \sqrt{x}) \nonumber \\
													&= {}_3F_4 \Big( \begin{array}{c} 1, \quad \frac\nu2+1,\quad \frac\nu2+\frac32 \\
																							 \frac32, \, \nu + \frac32,\, \frac{\nu-n}2+1, \, \frac{\nu-n}2+\frac32 \end{array}\, ; x\Big) \, \end{align*}
We recognize the coefficients in the associated Jensen polynomial \cite[p. 113]{Csordas}
   \[ \gamma_m = \frac{(-1)^m (1)_m (\frac\nu2+1)_m (\frac\nu2+\frac32)_m}
			           {(\frac32)_m (\nu+\frac32)_m (\frac{\nu-n}2+1)_m (\frac{\nu-n}2+\frac32)_m}	\,. \]
The related Jensen polynomial becomes the Laguerre-type hypergeometric polynomial
   \begin{align} \label{H4}
	    \mathcal{P}_s^{\mathbf H}(x;n) = \sum_{m=0}^s (-1)^m \binom sm \gamma_m\, x^m
			            = {}_4F_4 \Big( \begin{array}{cccc} -s,& 1, & \frac\nu2+1, & \frac\nu2+\frac32 \\
																							 \frac32, & \nu + \frac32, & \frac{\nu-n}2+1, & \frac{\nu-n}2+\frac32
																							 \end{array}\, ; x\Big) \,.
	 \end{align}
The special case $n = 0$ simplifies into
   \[ \mathcal{P}_s^{\mathbf H}(x;0) = {}_2F_2\left(-s, 1;\frac32, \nu + \frac32; x\right) \, , \]
while, for $n=1$ we have
	 \[ \mathcal{P}_s^{\mathbf H}(x;1) = {}_3F_3\left(-s,1,\frac\nu2+\frac32; \frac32,\nu + \frac32, \frac\nu2+\frac12; x\right) \,. \]
Now, by using the fact that for $\nu\in\left[-\frac{1}{2},\frac{1}{2}\right]$ the zeros of $\mathbf{H}_{\nu}$ and
$\mathbf{H}_{\nu}'$ are real, and also part {\bf a} of Theorem \ref{th3} it follows that the function $x \mapsto \mathcal H_{\nu, n}(2\sqrt{x})$ belongs to the
Laguerre-P\'olya class $\mathcal{LP}$ under assumption that $\nu\in\left[-\frac{1}{2},\frac{1}{2}\right]$ and $n\in\{0,1\},$ or $\nu\in\left(0,\frac{1}{2}\right]$ and $n=2.$ Consequently by using the well-known theorem of Jensen it follows that
the Jensen polynomial $\mathcal P_s^{\mathbf H}(x;n)$ has only real zeros. Now, the hypergeometric shape of the Jensen polynomial
$\mathcal P_s^{\mathbf H}(x;n)$ in the question is described in \eqref{H4}. Moreover, according to Csordas and Williamson \cite{csordas} the zeros of the Jensen polynomials are simple, and this completes the proof of the theorem.
\end{proof}

\begin{proof}[\bf Proof of Theorem \ref{th5}]
By equating the infinite series and the infinite product representation for $\mathbf{H}_{\nu}(z)$ which are given by (see \cite{bps})
$$\sqrt{\pi}2^{\nu}x^{-\nu-1}\Gamma\left(\nu+\frac{3}{2}\right)\mathbf{H}_{\nu}(x)=\prod_{n\geq 1}\left(1-\frac{x^2}{h_{\nu,n}^2}\right)$$
and
$$\mathbf{H}_{\nu}(x)=\left(\frac{x}{2}\right)^{\nu+1}\sum_{n\geq0}\frac{(-1)^n\left(\frac{x}{2}\right)^{2n}}
{\Gamma\left(n+\frac{3}{2}\right)\Gamma\left(n+\nu+\frac{3}{2}\right)},$$
we obtain
$${1\over{\Gamma({3\over 2})\Gamma(\nu+{3\over 2})}}-{{({x\over 2})^{2}}\over{{3\over 2}\Gamma({3\over 2})\Gamma(\nu+{3\over 2})}(\nu+{3\over 2})}+{\dots}
={2\over {\sqrt{\pi}\Gamma(\nu+{3\over 2})}}\left(1-{{x^{2}}\over{h_{\nu,1}^{2}}}\right)\left(1-{{x^{2}}\over{h_{\nu,2}^{2}}}\right){\dots}$$
or equivalently
$$1-\left({x\over 2}\right)^{2}{{2^{2}}\over{3(2\nu+3)}}+\left({x\over 2}\right)^{4}{{2^{4}}\over{3\cdot 5(2\nu+3)(2\nu+5)}}+{\dots}
=1-x^{2}\sum_{n\geq 1}{1\over {h_{\nu,n}^{2}}}+x^{4}\sum_{n\geq 1}{1\over {h_{\nu,n}^{2}}}\sum_{k\geq 1,\ k\ne n}{1\over {h_{\nu,k}^{2}}}+{\dots}.$$
The coefficients of the same power of $x$ must be equal, so the equality of the coefficients of $x^{2}$ proves the first relation in this theorem, and the equality of the coefficients of $x^{4}$ gives
$$\sum_{n\geq 1}{1\over {h_{\nu,n}^{2}}}\sum_{k\geq 1,\ k\ne n}{1\over {h_{\nu,k}^{2}}}={{1}\over{15(2\nu+3)(2\nu+5)}}.$$
Since the zeros $h_{\nu,n}$ are symmetric around the origin and because of the first Rayleigh sum of $h_{\nu,n}$-s, the previous equation becomes
$${1\over 2}\sum_{n\geq 1}{1\over {h_{\nu,n}^{2}}}\left(\sum_{k\geq 1}{1\over {h_{\nu,k}^{2}}}-{1\over{h_{\nu,n}^{2}}}\right)={{1}\over{15(2\nu+3)(2\nu+5)}}$$
and using again the first Rayleigh sum, it finally becomes the second Rayleigh sum as it is in the statement of the theorem.

To prove the relations on the Rayleigh sums of the zeros of $\mathbf{H}_{\nu}'$ we first equate the infinite sum and the infinite product representations of $\mathbf{H}_{\nu}'$ given by \eqref{sumS} and \eqref{prodS}, that is,
$$\sum_{n\geq0} {{(-1)^{n}(2n+\nu+1)({x\over 2})^{2n}}\over{(\nu+1)({3\over 2})_{n}(\nu+{3\over 2})_{n}}}=\prod_{n\geq 1}\left(1-{{x^{2}}\over{(h^{\prime}_{\nu,n})^{2}}}\right)$$
or equivalently
$$1-\left({x\over 2}\right)^{2}{{(\nu+3)}\over{(\nu+1)({3\over 2})_{1}({{\nu+3}\over 2})_{1}}}+\left({x\over 2}\right)^{4}{{(\nu+5)}\over{(\nu+1)({3\over 2})_{2}({{\nu+3}\over 2})_{2}}}+{\dots}$$ $$=1-x^{2}\sum_{n\geq 1}{1\over {(h_{\nu,n}^{\prime})^{2}}}+x^{4}\sum_{n\geq 1}{1\over {(h_{\nu,n}^{\prime})^{2}}}\sum_{k\geq 1,\ k\ne n}{1\over {(h_{\nu,k}^{\prime})^{2}}}+{\dots}.$$
The equality of the coefficients of the same power of $x$ on both sides gives the desired relations for the zeros of $\mathbf{H}_{\nu}'.$

Finally, to deduce the first two Rayleigh sums for the zeros of $\mathbf{H}_{\nu}''$ we proceed similarly as above.
\end{proof}

{\bf Acknowledgements}. The research of \'A. Baricz was supported by a grant of the Babe\c{s}-Bolyai University supporting excellence in scientific research, project number GSCE\underline{ }30246/2015.

\end{document}